\documentclass[10pt]{amsart}

\usepackage[T1,T2A]{fontenc} 
\usepackage[utf8x]{inputenc}
\usepackage[english]{babel} 
\usepackage{amsfonts, amsmath, amssymb, amsthm}
\usepackage{mathtools}
\usepackage{bm}
\usepackage{framed}
\usepackage{hyperref}
\hypersetup{
    colorlinks=true,
    citecolor=red,
    linkcolor=blue
  }

\usepackage[normalem]{ulem}

\allowdisplaybreaks[4]

\newtheorem{theorem}{Theorem}
\newtheorem{lemma}{Lemma}

\newtheorem{remark}{Remark}

\newtheorem{corollary}{Corollary}

\DeclarePairedDelimiter\floor{\lfloor}{\rfloor}

\DeclareMathOperator\corr{corr}
\DeclareMathOperator\var{var}

\begin{document}


\title[Parisian ruin with power-asymmetric variance near the
optimal point with application to\dots]{Parisian ruin with 
power-asymmetric variance near the optimal point with
application to many-inputs proportional reinsurance}

\author{Pavel Ievlev}
\address{
  Pavel Ievlev,
  Department of Actuarial Science,
  University of Lausanne\\
  UNIL-Dorigny, 1015 Lausanne, Switzerland
}
\email{ievlev.pn@gmail.com}

\bigskip

\date{\today}
\maketitle


{\bf Abstract:}
This paper investigates the Parisian ruin probability for
processes with power-asymmetric behavior of the variance near
the unique optimal point. We derive the exact asymptotics as
the ruin boundary tends to infinity and extend the previous
result \cite{DebickiEtAl2016} to the case when the length of
Parisian interval is of Pickands scale.
As a primary application, we extend the recent result 
\cite{Kepczynski2020} on the many inputs proportional 
reinsurance fractional Brownian motion risk model
to the Parisian ruin.







\section{Introduction}


  Consider the following reinsurance scheme: \( d \) companies
  share premiums and one claim process proportionally. Suppose
  that the risk process \( \bm{R} ( t ) \) is composed of a large
  number of i.i.d. sub-risk processes \( \bm{R}^{(i)} ( t ) \)
  representing independent businesses, and let each
  \( \bm{R}^{(i)} \) be driven by a fractional Brownian motion
  \( B_H ( t ) \). That is, let
  \[
    \bm{R}_N ( t )
    = \sum_{i = 1}^{N} \bm{R}^{(i)} ( t ),
    \quad \text{where} \quad
    \bm{R}^{(i)} ( t )
    = \bm{\alpha} + \bm{\mu} t - \bm{\sigma} B_H^{(i)} ( t ),
  \]
  where \( B_H^{(i)}, i \geq 1 \) are independent fractional
  Brownian motions and
  \( \bm{\alpha}, \bm{\mu}, \bm{\sigma} \in \mathbb{R}^d \).
  %
  %
  %
  In a recent contribution \cite{Kepczynski2020}, the
  authors derived the exact asymptotics of the simultaneous
  ruin probability
  \[
    \mathbb{P} \left\{
      \exists t \in [0, T] \colon
      \bm{R}_N ( t ) < 0
    \right\},
    \quad N \to \infty
  \]
  in the case of \( d = 2 \). In the present work we shall concentrate
  on the simultaneous Parisian ruin probability
  \begin{equation} \label{K-case}
    \Pi ( N )
    = \mathbb{P} \left\{
      \mathcal{P}_{[0, S], [0, T_N]} (\bm{R}_N) < 0
    \right\},
    \quad \text{where} \quad
    \mathcal{P}_{E, F} ( \bm{Z} )
    = \max_{i = 1, \dots, d}
    \sup_{t \in E}
    \inf_{s \in F}
    Z_i ( t + s ).
  \end{equation}



  Parisian stopping times have been first introduced in
  relation to barrier options in mathematical finance, see
  \cite{ChesneyPicquYor1997}, and since then attracted
  substantial interest.
  For the applications to actuarial risk theory, we refer to
  \cite{DassiosShanle2011semimarkov}, where risk process
  is treated as a surplus process of an insurance company
  with initial capital \( u \).
  
  In opposition to the well-studied classical ruin, when 
  the failure is recognized at the moment of surplus hitting 
  zero, the Parisian ruin is recognized only if the surplus 
  process has spent a sufficient, pre-specified amount of time 
  below zero.
  We refer to
  \cite{
    DassiosShanle2011, 
    CzrnaPalmowski2011,
    LoeffenCzarnaPalmowski2011, 
    CzrnaPalmowski2014-2, 
    CzarnaPalmowski2014}
  and the references therein for analysis of Parisian ruin in
  the one-dimensional L{\'e}vy surplus model.

  In the univariate Gaussian setup, Parisian ruin has been 
  investigated in \cite{DebickiEtAl2014} for self-similar
  Gaussian processes and in \cite{DebickiEtAl2016} for general
  Gaussian processes, satisfying some standard assumptions
  (see \cite{PiterbargBook}). Another interesting univariate
  case is of Parisian ruin over \textit{discrete sets}.
  In \cite{Jasnovidov2020}, the authors have proved that for
  the Brownian motion and \textit{equidistant grid} the
  asymptotics differs from the continuous one
  by some constant factor.
  
  There are many possible extensions of the notion of the
  Parisian ruin to multivariate risk processes, such as
  \textit{simultaneous Parisian ruin}, when all the components
  of a multivariate process should plunge below zero 
  \textit{at the same time} and remain there long enough 
  for the ruin to be attested. 
  This problem has recently been studied in \cite{Kriukov2021} 
  for the case when the risk process consists of two correlated 
  Brownian motions.
  Another possible extension is the \textit{joint} or 
  \textit{non-simultaneous Parisian ruin}, when ruin is
  attested if all the components experience Parisian ruin
  during some interval of time, but \textit{not necessarily at
  the same time}.
  This problem has been studied in \cite{Krystecki2021parisian},
  also for the bivariate Brownian motion with 
  \( \rho \in (-1,1)\).
  A third possible extension would be the notion of ``at least
  one'' ruin, suggested in the classical ruin context in
  \cite{Kepczynski2020}: the ruin is declared if
  \textit{either} of the processes has a Parisian ruin over
  time.



  In this paper, we derive the exact asymptotics as
  \( u \to \infty \)
  of the one-dimensional Parisian ruin probability
  \[
    \Pi ( u ) 
    = 
    \mathbb{P} \left\{ 
      \mathcal{P}_{[-S, S], [0, T_u]} \left( Z \right)
      > u
    \right\}
  \]
  for 
  \( T_u \to 0 \)
  at some specified rate and a class of Gaussian 
  processes with correlation structure
  \[
    \corr (Z(t), Z(s)) 
    =
    1 - D | t - s |^\alpha
    + o \left( | t - s |^\alpha \right)
  \]
  and a unique optimal point \( t_* = 0 \) of the variance
  with asymmetric behaviour near this point:
  \[
    \sigma ( t ) 
    = 
    1 
    - A_\pm | t |^{\gamma_\pm}
    \Big( 1 + o(1) \Big)
    \quad
    \text{as} \ t \to \pm 0,
  \]
  from which we further derive the exact asymptotics
  of the many-inputs Parisian ruin probability \eqref{K-case}.

  The asymptotic behaviour of 
  \( \Pi ( u ) \)
  for such class of processes is of interest by itself.
  Similar problems have recently been studied in 
  \cite{DebickiEtAl2014}
  and
  \cite{DebickiEtAl2016}.
  Our findings account for the previously discarded type of
  the
  Talagrand case with \( T > 0 \) (see Section
  \ref{main-results})
  and discover the new type of asymptotics therein:
  \[
    \Pi ( u ) 
    =
    e^{
      - \min \{ A_- T^{\gamma_-}, A_+ T^{\gamma_+} \}
    }
    \Psi ( u ),
  \]
  which rather surprisingly happens only if 
  \( \gamma_+ < \gamma_-, \alpha \)
  and does not happen if
  \( \gamma_- < \gamma_+, \alpha \).
    

  The paper is organized as follows. In Section
  \ref{main-results} we present our main findings.
  Theorem \ref{main-theorem} provides the exact asymptotics
  of 
  \( \Pi(u) \)
  for the general Gaussian process with power-asymmetric
  behaviour of the variance near the optimal point and
  under some assumption on the speed
  \( T_u \to 0 \)
  convergence. 
  It covers the previously unaccounted for case
  when the size of Parisian interval is equivalent to the
  Pickands scale of the process.
  Corollary \ref{MIPR-corollary} contains the exact
  asymptotics of the many-inputs Parisian ruin probability.
  The proof of Theorem \ref{main-theorem} is presented in a
  separate Section \ref{main-theorem-proof}.
  All known results and technical details are relegated to
  the Appendix.



\section{Main results} \label{main-results}

In this section, we first explain how to rewrite many-inputs
ruin probability in a form suitable for applying
Theorem \ref{main-theorem}, then specify the assumptions under
which the general theorem works and conclude with deriving the exact
asymptotics of the many-inputs proportional reinsurance ruin 
probability \ref{MIPR-corollary}.

Observe that by properties of Gaussian distribution
\[
  \bm{R}_N \stackrel{d}{=}
  \bm{\alpha} N + \bm{\mu} N t
  - \bm{\sigma} \sqrt{N} B_H ( t ),
\]
we can rewrite the ruin probability \eqref{K-case} as
\[
  \Pi_{[0, S]} ( N )
  = \mathbb{P} \left\{
    \mathcal{P}_{[0, S], [0, T_N]} ( Z )
    > \sqrt{N}
  \right\},
  \quad \text{where} \quad
  Z ( t ) = \frac{B_H ( t )}{D(t)},
  \quad D(t)
  = \max_{i = 1, \dots, d} ( \alpha_i + \mu_i t ).
\]


Next, we state a result concerning Parisian
ruin probabilities
\[
  \Pi ( u )
  = \mathbb{P} \left\{
    \mathcal{P}_{E, [0, T_u]} ( Z )
    > u
  \right\}
\]
for some large class of Gaussian processes
and then apply it to \( \Pi_{[0,S]} ( N ) \),
rewritten in the latter form.


\subsection{Assumptions}

Let \( E \) be a compact subset of \( \mathbb{R} \),
containing point \( 0 \) in its interior, and let
\( Z ( t ), t \in E \) be a centered Gaussian process with
a.s. continuous sample paths satisfying the following
two assumptions:
\begin{description}
  \item [Assumption A1] The standard variance function
    \( \sigma_Z \) of the Gaussian process \( Z \) attains its
    maximum on \( E \) at the unique point
    \( \widehat{\tau} = 0 \). Further, there exist positive
    constants \( \gamma_\pm \) and \( A_\pm \) such that
    \begin{equation}
      \label{eq:sigma-expansion}
      \sigma_Z ( t )
      =
      1 - A_\pm | t |^{\gamma_\pm}
      + o( | t |^{\gamma_\pm} )
      \quad \text{as} \quad
      t \to \pm 0.
    \end{equation}
  \item[Assumption A2] There exists some positive constant
    \( \alpha \in (0, 2] \)
    such that
    \[
      \corr \left(
	Z(t), Z(s)
      \right)
      = 1 - D \left| t - s \right|^{\alpha}
      + o \left( | t - s |^{\alpha} \right)
      \quad \text{as } t, s \to 0.
    \]
\end{description}

\begin{remark}
  Note that it follows from \textbf{A2} that there exists such
  \( \delta > 0 \) that
  \[
    \mathbb{E} \left\{ 
      \Big( 
	\overline{Z} ( t )
	-
	\overline{Z} ( s )
      \Big)^2
    \right\}
    <
    C | t - s |^\alpha
  \]
  for all \( t, s < \delta \).
\end{remark}


As it turns out, there are two numbers
\[
  \nu = \min \{ \alpha, \gamma_-, \gamma_+ \}
  \quad \text{and} \quad
  \gamma = \max \{ \gamma_-, \gamma_+ \}
\]
which determine \textit{the type} of the asymptotics, but before
proceeding to that, 
we also need the following assumption on the convergence
rate of
\(
  T_u \to 0
\):
\begin{description}
  \item [Assumption B]
    \(
      T_u = T u^{-2/\nu}
    \)
    for some \( T \in [0, \infty) \).
\end{description}

Next, we introduce two well-known and important constants in
the theory of Gaussian extremes, see
\cite{
  DebickiEtAl2017,
  DebickiEtAl2014,
  DebickiEtAl2016}.
Define for 
\( T \geq 0 \) 
and
\( \alpha \in (0, 2] \)
the generalized Pickands and Piterbarg constants
\[
  \mathcal{H}_\alpha^{\mathcal{P}} ( T )
  = 
  \lim_{\lambda \to \infty} 
  \frac{
    \mathcal{H}_{\alpha, 0}^{\mathcal{P}} ( \lambda, T )
  }{
    \lambda
  }
  \quad \text{and} \quad
  \mathcal{H}_{\alpha, h}^{\mathcal{P}} ( T )
  = \lim_{\lambda \to \infty}
  \mathcal{H}_{\alpha, h}^{\mathcal{P}} ( T ),
\]
where
\[
  \mathcal{H}_{\alpha, h}^{\mathcal{P}} ( T )
  =
  \mathbb{E}
  \exp \left( 
    \sup_{t \in [-\lambda, \lambda]}
    \inf_{s \in [0, T]}
    \Big( 
      \sqrt{2} B_{\alpha/2} ( t + s )
      - | t + s |^\alpha
      - h ( t + s )
    \Big)
  \right)
\]
for such continuous \( h \) that the limit exists.
%
We are in a position to formulate our main theorem.

\begin{theorem} \label{main-theorem}
  Let
  \(
    ( Z(t) )_{t \geq 0}
  \)
  be a centered Gaussian process satisfying assumptions
  \textbf{A1} and \textbf{A2}, and let \( T_u \) be a
  positive measurable function of \( u \) satisfying
  assumption \textbf{(B)}. Then
  \begin{description}
    \item [In the Pickands case
	\( \nu = \alpha \neq \gamma\)
      ] we have
      \[
	\Pi ( u )
	=
	C_S
	\ \mathcal{H}_\alpha^{\mathcal{P}}
	( D^{1/\alpha} T )
	\, u^{2/\nu - 2/\gamma} \,
	\Psi ( u )
	\Big( 1 + o(1) \Big),
      \]
      with
      \[
	C_S =
	A_+^{-1/\gamma_+} D^{1/\alpha}
	\Gamma \left( \frac{1}{\gamma_+} + 1 \right)
	1_{\gamma = \gamma_+}
	+ A_-^{-1/\gamma_-} D^{1/\alpha}
	\Gamma \left( \frac{1}{\gamma_-} + 1 \right)
	1_{\gamma = \gamma_-}.
      \]
    \item[In the Piterbarg case \( \nu = \alpha = \gamma \)] we have
      \[
	\Pi ( u )
	=
	\mathcal{H}_{\alpha, h}^{\mathcal{P}} ( D^{1/\alpha} T )
	\Psi ( u )
	(1 + o(1)),
      \]
      where
      \(
	h ( t )
	=
	A_- D^{-1/\alpha} | t |^{\gamma_-} 1_{t \leq 0}
	+ A_+ D^{-1/\alpha} | t |^{\gamma_+} 1_{t \geq 0}
      \).
    \item[In the Talagrand-1 case \( \gamma = \nu \neq \alpha \)]
      \[
	\Pi ( u )
	= C \Psi ( u )
	\Big( 1 + o ( 1 ) \Big),
	\quad
	C =
	\begin{cases}
	  1, & \gamma_+ \geq \gamma_-, \\
	  \exp(- \min \{ A_- T^{\gamma_-}, A_+ T^{\gamma_+}
	  \}), & \gamma_+ < \gamma_-.
	\end{cases}
      \]
  \end{description}
\end{theorem}



Now we proceed with our initial problem, to which end we
first have to study the behaviour of
\( \var B_H ( t ) / D ( t ) \).
Note that the derivative
\( \sigma' \)
of the variance function
\[
  \sigma ( t )
  =
  \frac{t^H}{D ( t )},
  \quad
  D ( t )
  = 
  \max_{i = 1, \dots, d} ( \alpha_i + \mu_i t )
\]
changes its sign exactly once, since
\[
  \sigma' ( t )
  = 
  \frac{t^{H-1}}{D^2 ( t )} G(t),
  \quad
  G(t)
  = 
  H D ( t ) - t D' (t),
\]
where
\( G(t) \)
is monotone and decreasing,
\( G(0) > 0 \)
and
\( G(t) \to - \infty \) as \( t \to \infty \).
Since
\( \sigma' \)
must change sign (possibly in a discontinuous manner) at the
optimal point 
\( t_* \) of \( \sigma \),
we have thus proved that such point is unique. Let us assume
that 
\( t_* \in (0, S) \) or \( S = \infty \),
since to account for the boundary maxima case an approach
slightly different to ours is needed.

The maximum can either be caused by intersection of some
two lines
\( l_\pm ( t ) = \alpha_\pm + \mu_\pm t \)
from
\( D(t) \), 
that is,
\[
  t_*
  = 
  \frac{\alpha_+ - \alpha_-}{\mu_- - \mu_+}
\]
in which case
\( \sigma' \) is discontinuous at \( t_* \),
or by a point away from the lines' intersections, satisfying
\( \sigma' ( t_* ) = 0 \),
that is,
\[
  t_*
  = 
  \frac{H \alpha}{\mu (1 - H)}.
\]
Finally, these two types of maxima can coincide, giving rise
to a power-asymmetric behavior near \( t_* \)
\[
  \frac{
    \sigma ( t )
  }{
    \sigma ( t_* )
  }
  = 
  1 - A_\pm | t - t_* |^{\gamma_\pm}
  \Big( 1 + o(1) \Big)
\]
with 
\( \gamma_\pm \in \{ 1, 2 \} \).
Precisely, if
\[
  -C_\pm^{(1)}
  = 
  \frac{
    \sigma_\pm' ( t_* )
  }{
    \sigma ( t_* )
  }
  = 
  \frac{H}{t_*}
  - \frac{\mu_\pm}{\alpha_\pm + \mu_\pm t_*}
  < 0,
\]
then 
\( \gamma_\pm = 1 \)
and
\( A_\pm = C_\pm^{(1)} \).
If on the other hand
\( C_\pm^{(1)} = 0 \),
then under the following non-degeneracy assumption
\[
  -C_\pm^{(2)} 
  = 
  \frac{
    \sigma_\pm'' ( t_* )
  }{
    \sigma ( t_* )
  }
  = 
  - \frac{H}{t_*^2}
  + \frac{\mu_\pm^2}{(\alpha_\pm + \mu_\pm t_*)^2}
  < 0
\]
we have
\( \gamma_\pm = 2 \)
and
\( A_\pm = C_\pm^{(2)} \).


Now we may introduce the natural asymptotic parameter
\[
  \widehat{N} = \frac{\sqrt{N}}{\sigma_Z ( t_* )}
\]
and formulate the corollary on the MIPR asymptotics.
\begin{corollary}[] \label{MIPR-corollary}
  Let \( T_N \) satisfy the condition
  \[
    \lim_{N \to \infty} T_N \widehat{N}^{1/H}
    = T \in [0, \infty).
  \]
  \begin{itemize}
    \item
      If either \( \gamma_+ \) or \( \gamma_- \) equals \( 2
      \),
      then
      \[
	\Pi ( N )
	=
	\frac{\sqrt{\pi}}{2}
	\frac{1}{\sqrt{A}}
	\Big(
	  1_{\gamma_+ = 2}
	  + 1_{\gamma_- = 2}
	\Big)
	\frac{
	  \mathcal{H}_{2H}^{\mathcal{P}} \left(
	    T / 2^{1/2H} t_*
	  \right)
	}{
	  2^{1/2H} t_*
	}
	\widehat{N}^{\zeta}
	\Psi ( \widehat{N} )
	\Big( 1 + o (1) \Big).
      \]
    \item
      If both \( \gamma_\pm = 1 > 2H \), then
      \[
	\Pi ( N )
	=
	\left(
	  \frac{1}{A_-}
	  + \frac{1}{A_+}
	\right)
	\frac{
	  \mathcal{H}_{2H}^{\mathcal{P}}
	  \left(
	    T / 2^{1/2H} t_*
	  \right)
	}{2^{1/2H} t_*}
	\widehat{N}^{\zeta}
	\Psi ( \widehat{N} )
	\Big( 1 + o(1) \Big).
      \]
    \item
      If both \( \gamma_\pm = 1 = 2H \), then
      \[
	\Pi ( N )
	= \mathcal{H}_{2H, h}^{\mathcal{P}} 
	( 
	  T / 2^{1/2H} t_* 
	)
	\Psi ( \widehat{N} )
	\Big( 1 + o(1) \Big).
      \]
    \item
      If \( \gamma_\pm = 1 < 2H \),
      then
      \[
	\Pi ( N )
	= \Psi ( \widehat{N} )
	\Big( 1 + o(1) \Big).
      \]
  \end{itemize}
\end{corollary}




\section{Proof of Theorem \ref{main-theorem}}
\label{main-theorem-proof}


This section is dedicated to the proof of Theorem
\ref{main-theorem}.


\subsection{Large vicinities.}
Looking ahead, we shall prove that only a small vicinity of
the optimal point contributes to the first order asymptotics,
and to evaluate its contribution we shall divide this small
vicinity into even smaller parts of some size \( q(u) \)
(referred to as the Pickands scale of the process \( Z \),
determined only by the covariance structure of \( Z \)), on
which the uniform local Pickands lemma may be applied.
It follows directly from the Piterbarg inequality and
the following obvious but important property of the Parisian
functional:
\begin{equation}
  \mathcal{P}_{E, F} ( f )
  = \sup_{t \in E} \inf_{s \in F} f ( t + s )
  \leq \sup_{t \in E} f ( t )
  \label{parisian-and-sup}
\end{equation}
that
\[
  \limsup_{u \to \infty}
  \frac{
    \Pi_{
      E \setminus [-\delta_-(u), \delta_+(u)]
    } ( u )
  }{
    u^\kappa \Psi ( u )
  }
  = 0
\]
for
\(
  \delta_\pm ( u )
  =
  u^{-2/\gamma_\pm}
  \ln^{2/\gamma_\pm} u
\)
and all \( \kappa > 0 \).
Since we intend to prove that
\(
  \Pi ( u )
  \sim \mathcal{H} u^\kappa \Psi ( u )
\)
for some \( \mathcal{H}, \kappa > 0 \), from this inequality
will follow that
\(
  \Pi ( u )
  \sim \Pi_{[-\delta_- ( u ), \delta_+ ( u )]} ( u )
\)
as \( u \to \infty \).
We can narrow the vicinity even further by once again using
\eqref{parisian-and-sup} and applying Lemma
\ref{large-vicinity-lemma} (Lemma 5.4 from
\cite{DebickiEtAl2019})
\[
  \limsup_{u \to \infty}
  \frac{
    \Pi_{
      [-\delta_-(u),\delta_+(u)]
      \setminus [
	-\Lambda u^{-2/\gamma_-},
	\Lambda u^{-2/\gamma_+}
      ]
    } ( u )
  }{
    u^\kappa \Psi ( u )
  }
  \leq C e^{-c \Lambda^\gamma}
\]
where
\( \gamma = \max \{ \gamma_-, \gamma_+ \} \). Due to this
inequality, we may concentrate on the exact asymptotics of
\[
  \Pi_{\Delta(u, \Lambda)} ( u )
  \quad \text{where} \quad
  \Delta ( u, \Lambda )
  = \Delta^+ ( u, \Lambda)
  \cup \Delta^- ( u, \Lambda ),
  \quad
  \Delta^\pm ( u, \Lambda )
  = \pm [0, \Lambda u^{-2/\gamma_\pm}]
\]
and then let \( \Lambda \to \infty \).

\subsection{Pickands intervals.}

Next, we introduce the left and right Pickands intervals
\( \Delta^\pm_k ( u, \lambda ) \)
with some additional parameter \( \lambda > 0 \)
\[
  \Delta^\pm_k ( u, \lambda )
  =
  \pm \lambda q(u) [ k , k + 1 ],
  \quad \text{where} \quad
  \nu = \min \{ \alpha, \gamma_+, \gamma_- \}, \
  q(u) = u^{-2/\nu},
\]
and the number of those fitting into the large vicinity
\( \Delta^\pm ( u, \Lambda ) \):
\[
  N^\pm ( u, \lambda, \Lambda ) = \floor*{
    \frac{
      | \Delta^\pm ( u ) |
    }{
      \lambda u^{-2/\nu}
    }
  }
  = \floor*{
    \frac{ \Lambda u^{\zeta_\pm} }{ \lambda }
  },
  \quad
  \zeta_\pm = \frac{2}{\nu} - \frac{2}{\gamma_{\pm}}
  = \max \left\{
    \frac{2}{\alpha} - \frac{2}{\gamma_\pm},
    \frac{2}{\gamma_\mp} - \frac{2}{\gamma_\pm},
    0
  \right\} \geq 0.
\]
We shall split the proof in four cases:
\begin{description}
  \item [Pickands case] \( \gamma \neq \nu = \alpha \)
  \item [Piterbarg case] \( \gamma = \nu = \alpha \)
  \item [Talagrand-1 case] \( \gamma = \nu \neq \alpha \)
  \item [Talagrand-2 case] \( \gamma \neq \nu \neq \alpha \)
\end{description}
In the Pickands case at least one of \( \zeta_\pm \) is nonzero,
therefore \( N^\pm \) grows as \( u^{\zeta_\pm} \).
In both Piterbarg and Talagrand-1 cases
\( \zeta_+ = \zeta_- = 0 \), hence
\( u \mapsto N^\pm \) is constant and we can set
\( \lambda = \Lambda \), in which case
\( N^\pm = 1 \) -- the zeroth Pickands interval coincides
with the informative vicinity.
The Talagrand-2 case is to be treated separately.

\subsection{Pickands case.}

To deal with the Pickands case, we employ the so-called double
sum method, which is based on the Bonferroni inequality
\[
  \bm{\Sigma}_1 ( u, \lambda, \Lambda )
  - \bm{\Sigma}_2 ( u, \lambda, \Lambda )
  \leq \Pi_{\Delta ( u, \Lambda )} ( u )
  \leq
  \bm{\Sigma}_1' ( u, \lambda, \Lambda ),
\]
where
\[
  \bm{\Sigma}_1 ( u, \lambda, \Lambda )
  =
  \underbrace{
    \sum_{k = 1}^{N^+ ( u, \lambda, \Lambda )}
    \Pi_{\Delta_k^+ ( u, \lambda )} ( u )
  }_{
    =: \bm{\Sigma}_1^+ ( u, \lambda, \Lambda )
  }
  + \underbrace{
    \sum_{k = 1}^{N^- ( u, \lambda, \Lambda )}
    \Pi_{\Delta_k^- ( u, \lambda, \Lambda )} ( u )
  }_{
    =: \bm{\Sigma}_1^- ( u, \lambda, \Lambda )
  }
  + \underbrace{
    \Pi_{
      \Delta_0^+ ( u, \lambda )
      \cup \Delta_0^- ( u, \lambda )
    } ( u )
  }_{
    =: \bm{\Sigma}_0 ( u, \lambda )
  },
\]
and 
\( \bm{\Sigma}_1'^{\pm} \) 
and
\(
  \bm{\Sigma}_1'
  = \bm{\Sigma}_1'^+ 
  + \bm{\Sigma}_1'^- 
  + \bm{\Sigma}_0
\) 
denote the same 
\( \bm{\Sigma}_1 \) 
but with 
\( N^\pm + 1 \) 
instead of
\( N^\pm \) 
(so that the collection of Pickands intervals indeed cover 
  \( \Delta ( u, \Lambda)\)),
and finally
\[
  \bm{\Sigma}_2 ( u, \lambda, \Lambda )
  = \sum_{
    \substack{
      \kappa, \kappa' \in \{ +, - \}, \\
      1 \leq k \leq N^\kappa ( u, \lambda, \Lambda ), \\
      1 \leq k' \leq N^{\kappa'} ( u, \lambda, \Lambda ), \\
      (\kappa, k) \neq (\kappa', k')
    }
  }
  \mathbb{P} \left\{
    \mathcal{P}_{
      \Delta^{\kappa}_{k} ( u, \lambda, \Lambda ),
      [0, T_u]
    } ( Z )
    > u,
    \mathcal{P}_{
      \Delta^{\kappa'}_{k'} ( u, \lambda, \Lambda ),
      [0, T_u]
    } ( Z )
    > u
  \right\}.
\]
Since the Parisian functional is bounded by \( \sup \)
functional, we can reduce the double sum estimate to the
classical (\( \sup \)) case
\[
  \bm{\Sigma}_2 ( u, \lambda, \Lambda )
  \leq
  \sum_{
    \substack{
      \kappa, \kappa' \in \{ +, - \}, \\
      1 \leq k \leq N^{\kappa} ( u, \lambda, \Lambda ), \\
      1 \leq k' \leq N^{\kappa'} ( u, \lambda, \Lambda ), \\
      (\kappa, k) \neq (\kappa', k')
    }
  }
  \mathbb{P} \left\{
    \sup_{
      \Delta^{\kappa}_{k} ( u, \lambda, \Lambda ),
    } Z ( t )
    > u,
    \sup_{
      \Delta^{\kappa'}_{k'} ( u, \lambda, \Lambda ),
    } Z ( t )
    > u
  \right\},
\]
and using similar arguments to those in \cite{JiRobert2018} obtain
\begin{equation}
  \label{double-sum-assertion}
  \lim_{\lambda \to \infty}
  \lim_{\Lambda \to \infty}
  \limsup_{u \to \infty}
  \frac{
    \bm{\Sigma}_2 ( u, \lambda, \Lambda )
  }{
    u^k \Psi ( u )
  }
  = 0
  \quad \text{for all } k > 0.
\end{equation}

We shall prove that if there exist two constants
\( C_\pm > 0 \)
such that
\[
  \lim_{\lambda \to \infty}
  \lim_{\Lambda \to \infty}
  \lim_{u \to \infty}
  \frac{
    \bm{\Sigma}_1^{\pm} ( u, \lambda, \Lambda )
  }{
    u^{\zeta_\pm} \Psi ( u )
  }
  = C_\pm,
\]
which together with the double sum estimate above and
\begin{equation}
  \lim_{\lambda \to \infty}
  \lim_{u \to \infty}
  \frac{\bm{\Sigma}_0 ( u, \lambda )}{\Psi ( u )}
  = \mathcal{H}_0
  \in (0, \infty),
  \label{zeroth-interval}
\end{equation}
yields
\begin{equation*}
  \lim_{\Lambda \to \infty}
  \lim_{u \to \infty}
  \frac{
    \Pi_{\Delta ( u, \Lambda )} ( u )
  }{
    u^{\zeta} \Psi ( u )
  }
  = C,
\end{equation*}
where
\(
  \zeta
  = \max \{ \zeta_+, \zeta_- \}
  > 0
\)
and
\(
  C = C_+ 1_{\zeta = \zeta_+}
  + C_- 1_{\zeta = \zeta_-}
\).
Finally, we obtain
\[
  \Pi ( u )
  \sim C u^{\zeta} \Psi ( u ).
\]

\subsection{Piterbarg and Talagrand-1 cases.}

As noted before, in both Piterbarg and Talagrand-1 cases
\(
  u \mapsto
  N^\pm ( u, \lambda, \Lambda )
\)
is constant. We may set 
\( \lambda = \Lambda \), 
which makes this constant equal one, and therefore
\[
  \Pi_{\Delta ( u, \Lambda )}
  = \bm{\Sigma}_0 ( u, \Lambda ).
\]
By \eqref{zeroth-interval}, we have
\[
  \lim_{\Lambda \to \infty}
  \lim_{u \to \infty}
  \frac{
    \bm{\Sigma}_0 ( u, \Lambda )
  }{
    \Psi ( u )
  }
  = \mathcal{H}_0
\]
for some \( \mathcal{H}_0 > 0 \),
which together with Lemma \ref{large-vicinity-lemma} yeilds
\(
  \Pi ( u ) \sim \mathcal{H}_0 \Psi ( u )
\).
In the Talagrand case \( \nu \neq \alpha \) we shall see that
\( \mathcal{H}_0 = 1 \), and therefore
\( \Pi ( u ) \sim \Psi ( u ) \).

\subsection{Talagrand-2 case.}

In the Talagrand-2 case we shall directly (that is, without
appealing to Pickands intervals) prove that there exists
positive and finite limit
\[
  \lim_{u \to \infty}
  \frac{\Pi_{\Delta(\lambda, u)} ( u )}{\Psi ( u )}
  =
  \begin{cases}
    1, & \gamma_- < \gamma_+, \\
    \exp(- \min \{ A_- T^{\gamma_-}, A_+ T^{\gamma_+} \} ),
    & \gamma_+ < \gamma_-,
  \end{cases}
\]
which ceases to depend on \( \lambda \) as long as 
\( \lambda > T \).



\subsection{Asymptotics of \( \Sigma_0 ( u, \lambda ) \) and the
Piterbarg and Talagrand-1 cases}

Denote
\[
  q ( u ) = u^{-2/\nu},
  \quad
  \nu = \min \{ \alpha, \gamma_+, \gamma_- \}.
\]
Note that
\[
  u^{2/\alpha} q(u)
  \xrightarrow[u \to \infty]{}
  1_{\nu = \alpha}
  \quad \text{and} \quad
  u^{2/\gamma_\pm} q(u)
  \xrightarrow[u \to \infty]{}
  1_{\nu = \gamma_\pm}.
\]

To apply the uniform local Pickands Lemma \ref{uLP}, let us rewrite the probability
\(
  \bm{\Sigma}_0 ( u, \lambda )
\)
in terms of a standardized process as follows:
\[
  \bm{\Sigma}_0 ( u, \lambda )
  = \mathbb{P} \left\{
    \mathcal{P}_{
      \Delta_0^+ ( u, \lambda )
      \cup \Delta_0^- ( u, \lambda )
    } ( Z )
    > u
  \right\}
  = \mathbb{P} \left\{
    \mathcal{P}' ( \xi_{u, 0} )
    > u
  \right\},
\]
where we have defined
\(
  \mathcal{P}' =
  \mathcal{P}_{[-\lambda,\lambda],[0,1]}
\)
and the family
\(
  \{ \xi_{u, 0} \colon u > 0 \}
\)
of centered Gaussian processes by
\[
  \xi_{u, 0} ( t, s )
  = \frac{
    \overline{Z} (
      q (u) t + T_u s 
    )
  }{
    1 + g (
      q(u) t + T_u s
    )
  },
  \quad
  1 + g ( t )
  = \frac{1}{ \sigma_Z ( t ) },
\]
Note that in contrast to
\( \xi_{u, k} \), \( k > 0 \) (see below), the
process \( \xi_{u, 0} \) is defined for
\( t \in [-\lambda, \lambda] \),
not \( [0, \lambda] \). This is neither a coincidence, nor a technical
decision: the two adjacent intervals near
the optimal point cannot be treated separately as it will be
evident from the result.


By \eqref{eq:sigma-expansion} and the definition of
\( g \) we have
\[
  u^2 g ( q(u) t + T_u s )
  \to
  h ( t + T s ) =
  h^+ ( t + T s ) 1_{\nu = \gamma_+}
  + h^- ( t + T s ) 1_{\nu = \gamma_-}.
\]
where
\[
  h^\pm ( \mu )
  = A_\pm | \mu |^{\gamma_\pm}
  1_{ \pm \mu > 0}.
\]

By assumption \textbf{A2} we have
\[
  u^2 \mathbb{E} \left\{
    \left| 
    Z( q(u) t + T_u s )
    - Z( q(u) t' + T_u s' )
    \right|^{2}
  \right\}
  \to 2 D 1_{\nu = \alpha} \Big|
  ( t - t' )
  + T (s - s')
  \Big|^\alpha,
\]
which means that the condition \textbf{C2}
\[
  \eta ( t, s )
  = B_{\alpha/2} ( t + s ) 1_{\nu = \alpha},
  \quad
  (t, s) \in [-D^{1/\alpha} \lambda, D^{1/\alpha} \lambda]
  \times [0, D^{1/\alpha} T]
\]
for \( T \geq 0 \).
%
By the uniform local Pickands Lemma \ref{uLP} (condition (C3) is obviously satisfied)
we have
\[
  \lim_{u \to \infty}
  \frac{\bm{\Sigma}_0 ( u, \lambda )}{\Psi ( u )}
  =
  \mathcal{H}_{\eta, h}^{\mathcal{P}_0}
  \Big(
    [-D^{1/\alpha} \lambda, D^{1/\alpha} \lambda]
    \times [0, D^{1/\alpha} T]
  \Big),
\]
where
\(
  \mathcal{H}_{\eta, h}^{\mathcal{P}_0} ( E )
  = \mathbb{E} \left\{ e^{\mathcal{P} ( \eta^h )} \right\}
\),
\( 
  \mathcal{P}_0
  = \mathcal{P}_{[-\lambda, \lambda], [0,T]}
\)
and
\begin{multline} \label{S0-process}
  \eta^h ( t, s )
  =
  \Big(
    \sqrt{2} B_{\alpha/2} ( t + s )
    - | t + s |^\alpha
  \Big) 1_{\nu = \alpha}
  - \\[7pt]
  - A_- D^{-\gamma_-/\alpha}
  | t + s |^{\gamma_-}
  1_{t + s \leq 0, \nu = \gamma_-}
  - A_+ D^{-\gamma_+/\alpha}
  | t + s |^{\gamma_+}
  1_{t + s \geq 0, \nu = \gamma_+}.
\end{multline}

\subsubsection{Piterbarg case.}

To prove the main theorem in the Piterbarg case
\( \gamma = \alpha = \nu \), that is when
\( \eta^h \) has all the terms
\[
  \eta^h ( t, s )
  =
  \Big(
    \sqrt{2} B_{\alpha/2} ( t + s )
    - | t + s |^\alpha
  \Big)
  - A_- D^{-\gamma_-/\alpha}
  | t + s |^{\gamma_-}
  1_{t + s \leq 0}
  - A_+ D^{-\gamma_+/\alpha}
  | t + s |^{\gamma_+}
  1_{t + s \geq 0},
\]
it remains to apply the standard result on the existence of
Piterbarg constants to see that
\[
  \lim_{\lambda \to \infty}
  \lim_{u \to \infty}
  \frac{\bm{\Sigma}_0 ( u, \lambda )}{\Psi(u)}
  =
  \lim_{\lambda \to \infty}
  \mathcal{H}_{\eta, h}^{\mathcal{P}_0}
  \Big(
    [-D^{1/\alpha} \lambda, D^{1/\alpha} \lambda]
    \times [0, D^{1/\alpha} T]
  \Big)
  = \mathcal{H}_{\eta, h}^{\mathcal{P}_0} (D^{1/\alpha} T)
\]
exists and is finite. This ends the proof of the main
theorem in the Piterbarg case.

\subsubsection{Talagrand-1 case.}

In the Talagrand-1 case \( \gamma = \nu \neq \alpha \) and,
therefore, the random part disappears from \eqref{S0-process}, whereas
all non-random terms are present:
\[
  \eta^h ( t, s )
  = - A_+ D^{-\gamma_+/\alpha}
  | t + s |^{\gamma_+}
  1_{t + s \geq 0}
  - A_- D^{-\gamma_-/\alpha}
  | t + s |^{\gamma_-}
  1_{t + s \leq 0}
\]
It remains to calculate
\( \mathcal{P}_0 ( \eta^h ) \)
explicitly: 
if
\( \lambda > T \),
we have
\begin{align*}
  \mathcal{P}_0 \left( \eta^h \right)
  &= 
  \sup_{t \in [-\lambda, \lambda]}
  \inf_{s \in [0, T]}
  \Big( 
    - A_+
    | t + s |^{\gamma_+}
    1_{t + s \geq 0}
    - A_-
    | t + s |^{\gamma_-}
    1_{t + s \leq 0}
  \Big)
  = \\[7pt] &=
  \max \left\{ 
    \sup_{t \in [-\lambda, -T]}
    \inf_{s \in [0, T]}
    \left( 
      -A_- | t + s |^{\gamma_-}
    \right)
    ,
    \sup_{t \in [-T, \lambda]}
    \inf_{s \in [0, T]}
    \left( 
      -A_+ | t + s |^{\gamma_+}
    \right)
  \right\}
  = \\[7pt] &=
  \max \left\{ 
    \sup_{t \in [-\lambda, -T]}
    (- A_- |t|^{\gamma_-} ),
    \sup_{t \in [-T, \lambda]}
    \left( 
      -A_+ | t + T |^{\gamma_+}
    \right)
  \right\}
  = \\[7pt] &=
  \max \left\{ 
    -A_- T^{\gamma_-}, 0
  \right\}
  = 0.
\end{align*}

Therefore, by lemma above in the Talagrand case
\[
  \mathcal{H}_{\eta, h}^{\mathcal{P}_0} \Big(
    [-D^{1/\alpha} \lambda, D^{1/\alpha} \lambda]
    \times [0, D^{1/\alpha} T]
  \Big)
  = 1
\]
for all \( \lambda > T \). Thus, we have proved that
\[
  \lim_{\lambda \to \infty}
  \lim_{u \to \infty}
  \frac{\bm{\Sigma}_0 ( u, \lambda )}{\Psi ( u )}
  = 1.
\]
This ends the proof of the main theorem in the Talagrand
case.


\subsection{Talagrand-2 case.}

To apply the uniform local Pickands Lemma \ref{uLP}, let us rewrite the probability
\(
  \Pi_{\Delta(u, \Lambda)} ( u )
\)
in terms of a standardized process as follows.
First, observe that the trivial equality
\[
  t + s
  = (t + s) 1_{t + s \geq 0}
  + (t + s) 1_{t + s \leq 0},
  \quad
  (t, s) \in [
    -\lambda u^{-2/\gamma_-},
    \lambda u^{-2/\gamma_+}
  ]
  \times [ 0, T_u ]
\]
may be rewritten as
\[
  t + s
  = q ( u, t', s' )
  = q_+ (u, t', s')
  + q_- (u, t', s'),
  \quad
  (t', s') \in [-\lambda, \lambda] \times [0, 1].
\]
where
\[
  q_\pm (u, t', s' )
  = \Big( u^{-2/\gamma_\pm} t' + T_u s' \Big)
  1_{\pm (u^{-2/\gamma_\pm} t' + T_u s' ) \geq 0}.
\]
Using this reparametrization, we rewrite
\[
  \Pi_{\Delta(u, \lambda)} ( u )
  = \mathbb{P} \left\{
    \mathcal{P}_{
      [-\lambda u^{-2/\gamma_-}, \lambda u^{-2/\gamma_+}]
    } ( Z )
    > u
  \right\}
  = \mathbb{P} \left\{
    \mathcal{P}' ( \xi_{u, 0} )
    > u
  \right\},
\]
where we have defined
\(
  \mathcal{P}' =
  \mathcal{P}_{[-\lambda,\lambda],[0,1]}
\)
and the family
\(
  \{ \xi_{u, 0} \colon u > 0 \}
\)
of centered Gaussian processes by
\[
  \xi_{u, 0} ( t', s' )
  = \frac{
    \overline{Z} (
      q(u, t', s')
    )
  }{
    1 + g (
      q(u, t', s')
    )
  },
  \quad
  1 + g ( t )
  = \frac{1}{ \sigma_Z ( t ) }.
\]

Note that
\[
  u^{2/\gamma_\pm}
  q_\pm ( u, t', s' )
  \to
  \Big( t' + 1_{\nu = \gamma_\pm} T s' \Big)
  1_{\pm (t' + 1_{\nu = \gamma_\pm} T s') \geq 0}
\]
uniformly in \( (t', s') \).
By \eqref{eq:sigma-expansion} and the definition of \(
g \) we have
\begin{multline*}
  u^2 g ( q(u, t', s') )
  \sim A_+ | u^{2/\gamma_+} q_+ (u, t', s') |^{\gamma_+}
  + A_- | u^{2/\gamma_-} q_- ( u, t', s' ) |^{\gamma_-}
  \sim \\[7pt] \sim
  A_+ | t' + 1_{\nu = \gamma_+} T s' |^{\gamma_+}
  1_{t' + 1_{\nu = \gamma_+} T s' \geq 0}
  + A_- | t' + 1_{\nu = \gamma_-} T s' |^{\gamma_-}
  1_{t' + 1_{\nu = \gamma_-} T s' \geq 0}.
\end{multline*}
By assumption \textbf{A2} we have
\[
  u^2 \mathbb{E} \left\{
    \left| 
    Z( q(u, t'_1, s'_1) )
    - Z( q(u, t'_2, s'_2 ) )
    \right|^{2}
  \right\}
  \to
  0,
\]
which means that the condition \textbf{C2}
\(
  \eta ( t, s )
  = 0
\).

By the uniform local Pickands Lemma \ref{uLP} (condition
\textbf{C3} is obviously satisfied) we
have
\[
  \lim_{u \to \infty}
  \frac{\Pi_{\Delta(u, \lambda)} ( u )}{\Psi ( u )}
  = \mathcal{H}_{0, h}^{\mathcal{P}_0}
  ([-\lambda, \lambda] \times [0, T]),
\]
where
\(
  \mathcal{H}_{0,h}^{\mathcal{P}_0} ( E )
  = e^{\mathcal{P}_0 (-h)}
\),
\( 
  \mathcal{P}_0
  = \mathcal{P}_{[-\lambda, \lambda], [0, T]}
\)
and
\[
  h ( t, s )
  = A_- | t + 1_{\nu = \gamma_-} s |^{\gamma_-}
  1_{t + 1_{\nu = \gamma_-} s \leq 0}
  + A_+ | t + 1_{\nu = \gamma_+} s |^{\gamma_+}
  1_{t + 1_{\nu = \gamma_+} s \geq 0}.
\]

Since we are looking at a case where
\( \gamma \neq \nu \), 
either
\( 1_{\nu = \gamma_-} \) or \( 1_{\nu = \gamma_+} \)
is zero.
Suppose, \( \nu = \gamma_- \neq \gamma_+ \). Then
\[
  h ( t, s )
  =
  -A_-
  | t + s |^{\gamma_-}
  1_{t + s \leq 0}
  - A_+ | t |^{\gamma_+} 1_{t \geq 0}.
\]
therefore, we have
\begin{multline*}
  \mathcal{P} ( h )
  =
  \max \left\{
    \sup_{t \in [-\lambda, 0]}
    \inf_{s \in [0, t]}
    \big(
      -a_- | t + s |^{\gamma_-}
      1_{t + s \leq 0}
    \big),
    \sup_{t \in [0, \lambda]}
    \inf_{s \in [0, t]}
    \big(
      - a_+ | t |^{\gamma_+}
    \big)
  \right\}
  = \\[7pt] =
  \max \left\{
    \sup_{t \in [-\lambda, 0]}
    \inf_{\mu \in [t, t + t]}
    \big(
      - a_- | \mu |^{\gamma_-}
      1_{\mu \leq 0}
    \big),
    0
  \right\}
  = 0
\end{multline*}
since the first term is non-positive.

If, on the other hand, \( \nu = \gamma_+ \neq \gamma_- \),
we have
\[
  h ( t , s )
  =
  - a_- | t |^{\gamma_-} 1_{t \leq 0}
  - a_+ | t + s |^{\gamma_+} 1_{t + s \geq 0}
\]
and therefore for \( \lambda > T \) we have
\begin{align*}
  \mathcal{P}(h)
  &= \max \left\{
    \sup_{t \in [-\lambda, 0]}
    \inf_{s \in [0, T]}
    \Big(
      -A_- | t |^{\gamma_-}
      -A_+ | t + s |^{\gamma_+}
      1_{t + s \geq 0}
    \Big),
    \sup_{t \in [0, \lambda]}
    \inf_{s \in [0, T]}
    \Big(
      -A_+ | t + s |^{\gamma_+}
    \Big)
  \right\}
  = \\[7pt]
  &=
  \max \left\{
    \sup_{t \in [-\lambda, 0]}
    \inf_{\mu \in [t, t + T]}
    \Big(
      -A_- | t |^{\gamma_-}
      -A_+ | \mu |^{\gamma_+} 1_{\mu \geq 0}
    \Big),
    - A_+ | T |^{\gamma_+}
  \right\}
  = \\[7pt]
  &=
  \max \left\{
    \sup_{t \in [-\lambda, 0]}
    \min \left\{
      -A_- | t |^{\gamma_-},
      \inf_{\mu \in [0, \max (t + T, 0)] }
      \Big(
	-A_- | t |^{\gamma_-}
	- A_+ | \mu |^{\gamma_+}
      \Big)
    \right\},
    -A_+ | T |^{\gamma_+}
  \right\}
  = \\[7pt]
  &=
  \max \left\{
    \sup_{t \in [-\lambda, 0]}
    \Big(
      -A_- | t |^{\gamma_-}
      - A_+ | \max ( t + T, 0) |^{\gamma_+}
    \Big),
    -A_+ | T |^{\gamma_+}
  \right\}
  = \\[7pt]
  &=
  \max \left\{
    \sup_{t \in [-\lambda, -T]}
    \Big(
      -A_- | t |^{\gamma_-}
    \Big),
    \sup_{t \in [-T, 0]}
    \Big(
      -A_- | t |^{\gamma_-}
      - A_+ | t + T |^{\gamma_+}
    \Big),
    -A_+ | T |^{\gamma_+}
  \right\}
  = \\[7pt]
  &=
  \max \left\{
    -A_- | T |^{\gamma_-},
    \sup_{\mu \in [0, T]}
    \Big(
      -A_- | \mu |^{\gamma_-}
      - A_+ | T - \mu |^{\gamma_+}
    \Big),
    -A_+ | T |^{\gamma_+}
  \right\}
  = \\[7pt]
  &=
  - \min \left\{
    A_- T^{\gamma_-},
    A_+ T^{\gamma_+}
  \right\}
\end{align*}

We have thus proved that
\[
  \mathcal{H}_{0, h}^{\mathcal{P}}
  ([-\lambda, \lambda] \times [0, T])
  = e^{\mathcal{P}(-h)}
  =
  \begin{cases}
    1, & \gamma_- < \gamma_+, \\
    \exp(- \min \{ A_- T^{\gamma_-}, A_+ T^{\gamma_+} \}
    ), & \gamma_+ < \gamma_-
  \end{cases}
\]
for all \( \lambda > T \).



\subsection{Pickands case.}

Now we proceed to the Pickands case.

To find the aforementioned asymptotics of
\(
  \bm{\Sigma}_1^\pm ( u, \lambda, \Lambda )
\)
we shall first find the uniform in
\(
  k \in \{ 1, \dots, N^\pm ( u, \lambda, \Lambda ) \}
\)
asymptotics of each summand
\(
  \Pi_{\Delta^\pm ( u, \lambda, \Lambda )} ( u )
\)
and then sum them up. To this end, let us rewrite the probability
in the form required for the uniform local Pickands Lemma \ref{uLP}
\[
  \Pi_{\Delta_k^\pm ( u, \lambda, \Lambda )} ( u )
  = \mathbb{P} \left\{
    \mathcal{P}_{
      \Delta_k^\pm ( u, \lambda, \Lambda ),
      [0, T_u]
    } ( Z )
    > u
  \right\}
  = \mathbb{P} \left\{
    \mathcal{P}'' ( \xi_{u, k}^\pm )
    > u
  \right\},
\]
where we have defined
\(
  \mathcal{P}''
  = \mathcal{P}_{[0, \lambda], [0, 1]}
\)
and the family
\[
  \{
    \xi_{u, k}^{\kappa} \colon
    \kappa \in \{ +, - \},
    1 \leq k \leq N^\pm ( u, \lambda, \Lambda ),
    u > 0
  \}
\]
of centered Gaussian processes by
\[
  \xi_{u, k}^\pm ( t, s ) = \frac{
    \widetilde{Z}_{u, k}^\pm ( t, s )
  }{1 + h_{u, k}^\pm ( t, s )},
  \quad
  1 + h_{u, k}^\pm ( t, s )
  =
  \frac{1}{
    \sigma_{Z}(
      \pm q ( u ) ( \lambda k + t )
      + T_u s
    )
  }
\]
and
\[
  \widetilde{Z}_{u, k}^\pm ( t, s )
  = \frac{
    Z (
      \pm q ( u ) ( \lambda k + t )
      + T_u s
  )}
  {
    \sigma_Z (
      \pm q ( u ) ( \lambda k + t )
      + T_u s
    )
  },
\]
\( t \in [0, \lambda], s \in [0, 1] \). Note that
\(
  \widetilde{Z}_{u, k}^\pm
\)
is a centered Gaussian random field with unit variance and
continuous paths. Besides,
\(
  h_{u, k}^\pm \in C_0 ( [0, \lambda] \times [0, T] )
\),
that is, it is a continuous function on
\( [0, \lambda] \times [0, T] \),
such that
\(
  h_{u, k}^\pm ( 0, 0 ) = 0
\).

There is, however, a pitfall in trying to apply uLP directly
to \( \xi_{u, k}^\pm \). Even though \( u^2 h_{u, k} \) may have a
limit \( h \) for each \( k \), it is never uniform. In
other words, the condition (C1)
\[
  \lim_{u \to \infty}
  \sup_{k \in K_u, (t, s) \in [0, \lambda] \times [0, T]}
  \left|
  u^2 h_{u, k} ( t, s ) - h ( t, s )
  \right|
  = 0
\]
is not satisfied.
To get around this inconvenience, we shall coarsen the inequality describing the
event by taking \( 1 + h_{u, k} \) out of \( \mathcal{P}'' \)
\[
  \mathbb{P} \left\{
    \mathcal{P}'' ( \widetilde{Z}_{u, k} )
    / U ( 1 + h_{u, k}^\pm )
    > u
  \right\}
  \leq
  \Pi_{\Delta_k^\pm ( u, \lambda )} ( u )
  =
  \mathbb{P} \left\{
    \mathcal{P}'' ( \xi_{u, k} )
    > u
  \right\}
  \leq
  \mathbb{P} \left\{
    \mathcal{P}'' ( \widetilde{Z}_{u, k} )
    / L ( 1 + h_{u, k}^\pm )
    > u
  \right\},
\]
where
\[
  U ( f )
  = \sup_{t \in [0, \lambda]}
  \sup_{s \in [0, T]}
  f ( t + s )
  \quad \text{and} \quad
  L ( f )
  = \inf_{t \in [0, \lambda]}
  \inf_{s \in [0, T]}
  f ( t + s ).
\]
Let us rewrite it in a handier fashion as
\[
  \mathbb{P} \left\{
    \mathcal{P}'' ( \widetilde{Z}_{u, k} )
    > u
    \Big( 1 + U ( h_{u, k}^\pm ) \Big)
  \right\}
  \leq
  \Pi_{\Delta_k^\pm ( u, \lambda )}
  = \mathbb{P} \left\{
    \mathcal{P}'' ( \xi_{u, k})
    > u
  \right\}
  \leq
  \mathbb{P} \left\{
    \mathcal{P}'' ( \widetilde{Z}_{u, k} )
    > u
    \Big( 1 + L ( h_{u, k}^\pm ) \Big)
  \right\}.
\]
Using the assumption that \( A_\pm > 0 \), we get
\[
  L ( h_{u, k}^+ )
  = A_+ \Big| 
  \lambda q ( u ) k
  \Big|^{\gamma_+},
  \quad U ( h_{u, k}^+ )
  = A_+ \Big|
  \lambda q ( u ) (k + 1) + T_u
  \Big|^{\gamma_+}
\]
and
\[
  L ( h_{u, k}^- )
  = A_- \Big| 
  - \lambda q(u) k + T_u
  \Big|^{\gamma_-},
  \quad
  U ( h_{u, k}^- )
  = A_- \Big| 
  \lambda q(u) (k + 1)
  \Big|^{\gamma_-}.
\]
All four bounds can be rewritten as follows:
\[
  L ( h_{u, k}^\pm )
  = A_\pm \Big| 
  \lambda q(u) k 
  - q(u) m_\pm (u)
  \Big|^{\gamma_\pm},
  \quad
  U ( h_{u, k}^\pm )
  = A_\pm \Big| 
  \lambda q(u) k
  + q(u) ( \lambda + m_\mp ( u ) )
  \Big|^{\gamma_\pm},
\]
where 
\( m_+ ( u ) = 0 \)
and
\( m_- ( u ) = T_u / q(u) \to T \)
as
\( u \to \infty \).
It is important for us that \( m_\pm \) do not depend on 
\( k \) and have finite limits as \( u \to \infty \).

In order to apply the uniform local Pickands Lemma \ref{uLP} to the upper bound, we set
\[
  g_{u, k} ( \lambda )
  = u
  \Big(
    1 + A_\pm 
    \Big| 
    \lambda q ( u ) k 
    - q (u) m_\pm ( u )
    \Big|^{\gamma_\pm}
  \Big)
\]
and note that the condition (C2) remains valid with
\( g_{u, k} \) instead of \( u \). It now follows
directly from the uniform local Pickands Lemma \ref{uLP} that with
\( 
  \mathcal{P}
  = \mathcal{P}_{[0, \lambda], [0, T]}
\)
we have
\[
  \mathbb{P} \left\{
    \mathcal{P}'' ( \widetilde{Z}_{u, k} )
    > g_{u, k} ( \lambda )
  \right\}
  = \mathcal{H}_{2H}^{\mathcal{P}}
  \left(
    D^{1/\alpha} \lambda,
    D^{1/\alpha} T
  \right)
  \Psi ( u )
  \exp \left(
    - A_\pm \Big|
    \lambda u^{-\zeta_\pm} k
    - u^{-\zeta_\pm} m_\pm ( u )
    \Big|^{\gamma_+}
  \right)
  (1 + o (1) ),
\]
where \( o(1) \) is uniform in
\( k \in \{ 1, \dots, N^\pm ( u, \lambda, \Lambda) \} \).
Thus,
\begin{multline*}
  \bm{\Sigma}_1^\pm ( u, \lambda, \Lambda )
  \leq
  \sum_{k = 1}^{N^\pm ( u, \lambda, \Lambda )}
  \mathbb{P} \left\{
    \mathcal{P} ( \widetilde{Z}_{u,k} )
    >
    g_{u, k} ( \lambda )
  \right\}
  \sim \\[7pt] \sim
  \mathcal{H}_{2H}^{\mathcal{P}}
  \left(
    D^{1/\alpha} \lambda,
    D^{1/\alpha} T
  \right)
  \Psi (u)
  \sum_{k = 1}^{N^\pm ( u, \lambda, \Lambda )}
  \exp \left(
    - A_\pm \Big|
    \lambda u^{-\zeta_\pm} k
    - u^{-\zeta_\pm} m_\pm ( u )
    \Big|^{\gamma_\pm}
  \right)
  \sim \\[7pt] \sim
  \frac{
    \mathcal{H}_{2H}^{\mathcal{P}}
    \left(
      D^{1/\alpha} \lambda,
      D^{1/\alpha} T
    \right)
  }{
    \lambda
  }
  u^{\zeta_\pm} \Psi ( u )
  \int_0^{\Lambda}
  e^{-A_\pm x^{\gamma_\pm}}
  \mathop{dx}
\end{multline*}
and, letting \( u \to \infty \), then \( \Lambda \to
\infty \) and finally \( \lambda \to \infty \), we see
that for large enough \( u \), holds
\[
  \lim_{\lambda \to \infty}
  \lim_{\Lambda \to \infty}
  \limsup_{u \to \infty}
  \frac{
    \bm{\Sigma}_1^\pm ( u, \lambda, \Lambda )
  }{
    u^{\zeta_\pm} \Psi ( u )
  }
  \leq
  \Gamma \left( \frac{1}{\gamma_\pm} + 1 \right)
  A_\pm^{-1/\gamma_\pm}
  D^{1/\alpha}
  \mathcal{H}_{\alpha}^{\mathcal{P}} ( D^{1/\alpha} T )
\]
where
\[
  \mathcal{H}_{\alpha}^{\mathcal{P}} ( T )
  = \lim_{\lambda \to \infty}
  \frac{
    \mathcal{H}_{\alpha}^{\mathcal{P}}
    ( \lambda, T )
  }{\lambda}
  \in (0, \infty).
\]
Using lower bound in much the same fashion, we obtain
\[
  \lim_{\lambda \to \infty}
  \lim_{\Lambda \to \infty}
  \limsup_{u \to \infty}
  \frac{
    \bm{\Sigma}_1^\pm ( u, \lambda, \Lambda )
  }{
    u^{\zeta_\pm} \Psi ( u )
  }
  =
  \Gamma \left( \frac{1}{\gamma_\pm} + 1 \right)
  A^{-1/\gamma_\pm} D^{1/\alpha}
  \mathcal{H}_{\alpha}^{\mathcal{P}} (D^{1/\alpha} T).
\]
The same formula obviously holds for 
\( \bm{\Sigma}_1'^\pm \). 
To conclude the proof in the Pickands case, it remains to notice that
\begin{multline*}
  \lim_{\lambda \to \infty}
  \lim_{\Lambda \to \infty}
  \lim_{u \to \infty}
  \frac{
    \bm{\Sigma}_1 ( u, \lambda, \Lambda )
  }{
    u^\zeta \Psi ( u )
  }
  = \\[7pt] =
  \lim_{\lambda \to \infty}
  \lim_{\Lambda \to \infty}
  \lim_{u \to \infty}
  \left(
    u^{\zeta_+ - \zeta}
    \frac{
      \bm{\Sigma}_1^+ ( u, \lambda, \Lambda )
    }{u^{\zeta_+} \Psi(u)}
    + u^{\zeta_- - \zeta}
    \frac{
      \bm{\Sigma}_1^- ( u, \lambda, \Lambda )
    }{u^{\zeta_-} \Psi(u)}
    + u^{-\zeta}
    \frac{
      \bm{\Sigma}_0 ( u, \lambda )
    }{\Psi(u)}
  \right)
  = \\[7pt] =
  \left(
    \Gamma \left( \frac{1}{\gamma_+} + 1 \right)
    A_+^{-1/\gamma_+}
    1_{\zeta = \zeta_+}
    + \Gamma \left( \frac{1}{\gamma_-} + 1 \right)
    A_-^{-1/\gamma_-}
    1_{\zeta = \zeta_+}
  \right)
  D^{1/\alpha}
  \mathcal{H}_\alpha^{\mathcal{P}} ( D^{1/\alpha} T ),
\end{multline*}
that the same is obviously true for \( \bm{\Sigma}_1' \),
and that 
\( 1_{\zeta = \zeta_\pm} = 1_{\gamma = \gamma_\pm} \).




\section{Appendix}


In this appendix, we recall some known results necessary for
the proofs of previous section.



\subsection{Parisian functional continuity}

Let us show that the Parisian functional
\( \mathcal{P} \colon C ( E \times F ) \to \mathbb{R} \)
is continuous in uniform topology. 
To this end, we take an arbitrary function
\( 
  f \in C ( E \times F )
\)
and a family
\( 
  \{ f_\varepsilon, \varepsilon > 0 \}
  \subset C ( E \times F )
\),
which converges to a function
\( 
  f \in C ( E \times F )
\)
uniformly
\[
  \sup_{(t, s) \in E \times F}
  \Big| 
    f ( t, s )
    - 
    f_\varepsilon ( t, s )
  \Big|
  < \varepsilon
\]
as 
\( \varepsilon \to 0 \),
from which we obtain
\[
  \begin{cases}
    f ( t, s ) - f_\varepsilon ( t, s ) < \varepsilon,
    \\
    f_\varepsilon ( t, s ) - f ( t, s ) < \varepsilon
  \end{cases}
  \quad
  \text{for all} 
  \quad
  (t, s) \in E \times F.
\]
Hence,
\[
  \sup_{t \in E} \inf_{s \in F} 
  f (t, s)
  < 
  \varepsilon 
  + \sup_{t \in E} \inf_{s \in F} 
  f_\varepsilon (t, s)
  \quad \text{and} \quad
  \sup_{t \in E} \inf_{s \in F} 
  f_\varepsilon (t, s)
  <
  \varepsilon
  + \sup_{t \in E} \inf_{s \in F} f (t, s),
\]
or, equivalently,
\[
  \Big| 
    \sup_{t \in E} \inf_{s \in F} f (t, s)
    -
    \sup_{t \in E} \inf_{s \in F} f_\varepsilon (t, s)
  \Big|
  =
  \Big| 
    \mathcal{P} ( f ) - \mathcal{P} ( f_\varepsilon )
  \Big|
  < \varepsilon.
\]



\subsection{Uniform local Pickands lemma}

The following lemma is from \cite{DebickiEtAl2017}, it is
reproduced here for the reader's convenience.
Let
\[
  \xi_{u, \tau_u} ( t )
  = \frac{Z_{u, \tau_u} ( t )}{1 + h_{u, \tau_u} ( t )},
  \quad t \in E,
  \tau_u \in K_u,
\]
be a family of centered Gaussian random fields with
\( Z_{u, \tau_u} \) a centered Gaussian random field with unit
variance and continuous paths, and \( h_{u, \tau_u} \)
belonging to \( C_0 ( E ) \), that is, \( h_{u, \tau_u} \) is
a continuous function on \( E \), such that \( h_{u, \tau_u} (
0 ) = 0 \). We assume that \( E \) is a compact subset of \(
\mathbb{R}^d \) and \( 0 \in E \).

The Parisian functional \( \Gamma_{E, F} \)
\[
  \Gamma_{E, F} ( X )
  = \sup_{t \in E} \inf_{s \in F} X ( t + s )
\]
satisfies the conditions
\begin{description}
  \item [(F1)] there exists \( c > 0 \) such that
    \(
      \Gamma ( f ) \leq c \sup_{t \in E} f ( t )
    \)
    for any \( f \in C(E) \)
  \item[(F2)] \( \Gamma ( a f + b ) = a \Gamma ( f ) + b \)
    for any \( f \in C(E) \) and \( a > 0 \), \( b \in
    \mathbb{R} \)
\end{description}
of the paper \cite{DebickiEtAl2017}. Therefore, under
conditions \( (C0) - (C3) \)
\begin{description}
  \item [(C0)] \( \lim_{u \to \infty} \inf_{\tau_u \in K_u}
    g_{u, \tau_u} = \infty \)
  \item[(C1)]  there exists \( h \in C_0(E) \) such that
    \[
      \lim_{u \to \infty} \sup_{\tau_u \in K_u, t \in E}
      \left|
      g_{u, \tau_u}^2 h_{u, \tau_u} ( t ) - h ( t )
      \right|
      = 0
    \]
  \item[(C2)] there exists \( \theta_{u, \tau_u} ( t, s ) \)
    such that
    \[
      \lim_{u \to \infty} \sup_{\tau_u \in K_u}
      \sup_{s \neq t \in E}
      \left|
      g_{u, \tau_u}^2
      \frac{
	\mathbb{E} \left\{ 
	  \left| 
	  Z_{u, \tau_u} ( t ) - Z_{u, \tau_u} ( s )
	  \right|^{2}
	\right\}
      }{2 \theta_{u, \tau_u} ( t, s )}
      - 1
      \right|
      = 0,
    \]
    where for some centered Gaussian random field
    \( \eta ( t ) \), \( t \in \mathbb{R}^d \)
    with continuous paths and \( \eta ( 0 ) = 0 \),
    \[
      \lim_{u \to \infty} \sup_{\tau_u \in K_u}
      \left|
      \theta_{u, \tau_u} ( t, s )
      - \mathbb{E} \left\{
	\left| 
	\eta ( t ) - \eta ( s )
	\right|^{2}
      \right\}
      \right|
      = 0.
    \]
  \item[(C3)] there exists \( a > 0 \) such that
    \[
      \limsup_{u \to \infty} \sup_{\tau_u \in K_u}
      \sup_{s \neq t \in E}
      \frac{\theta_{u, \tau_u} ( t, s )}{
	\sum_{i = 1}^d | s_i - t_i |^a
      }
      < \infty,
    \]
    \[
      \lim_{\varepsilon \to 0} \limsup_{u \to \infty}
      \sup_{\tau_u \in K_u}
      \sup_{\| t - s \| < \varepsilon, t, s \in E}
      g_{u, \tau_u}^2 \mathbb{E} \left\{ 
	\left[
	  Z_{u, \tau_u} ( t ) - Z_{u, \tau_u} ( s )
	\right] 
	Z_{u, \tau_u} ( 0 ) 
      \right\}
      = 0.
    \]
\end{description}
we have
\begin{theorem}[]
  \label{uLP}
  Under assumptions (C0)-(C3), if, further,
  \( \mathbb{P} ( \Gamma_{E, F} (\xi_{u, \tau_u}) >
  g_{u, \tau_u} ) > 0 \) for all \( \tau_u \in K_u \) and all
  large \( u \), then
  \[
    \lim_{u \to \infty} \sup_{\tau_u \in K_u}
    \left|
    \frac{
      \mathbb{P} \{
	\Gamma ( \xi_{u, \tau_u} )
	> g_{u, \tau_u}
      \}
    }{
      \Psi ( g_{u, \tau_u} )
    }
    - \mathcal{H}_{\eta, h}^\Gamma ( E )
    \right| = 0,
  \]
  where
  \[
    \mathcal{H}_{\eta, h}^G ( E )
    = \mathbb{E} \left\{ e^{\Gamma (\eta^h)} \right \},
    \quad
    \eta^h ( t )
    = \sqrt{2} \eta ( t )
    - \sigma^2_\eta ( t )
    - h ( t ).
  \]
\end{theorem}



\subsection{Large vicinity cut-off lemma}
\label{large-vicinity-lemma}

Next lemma is from \cite{DebickiEtAl2019} (Lemma 5.4), but
instead of the version therefrom, we give a
version suitable for our needs.
This lemma allows one to get rid of the complement of the 
Piterbarg vicinity in all three (Piterbarg, Pickands and
Talagrand) cases (see proof of Theorem \ref{main-theorem}
below).

\begin{lemma}[]
  There exist positive constants 
  \( C \), \( c \), \( u_0 \) and \( \Lambda_0 \)
  such that for
  \( \Lambda \geq \Lambda_0 \)
  and
  \( u \geq u_0 \)
  \[
    \mathbb{P} \left\{ 
      \exists t \in \left[ 
	-\delta_- ( u ), \Lambda u^{-2/\gamma_-}
      \right] 
      \cup \left[ 
	\Lambda u^{-2/\gamma_+}, \delta_+ ( u )
      \right]
      \colon
      Z ( t )
      > u
    \right\}
    \leq
    C e^{- c \Lambda^\gamma}
    \mathbb{P} \left\{ 
      Z ( 0 ) > u
    \right\}.
  \]
  where
  \( 
    \delta_\pm ( u ) 
    = u^{-2/\gamma_\pm} \ln^{2/\gamma_\pm} u.
  \)
\end{lemma}




\section{Acknowledgments}

Financial support by SNSF Grant 200021-196888
is kindly acknowledged.



\bibliographystyle{ieeetr}
\bibliography{$BIB}


\end{document}